\catcode`\@=11
\documentclass[11pt]{amsart}

\usepackage{sasha}
\usepackage{sasha_prih}
\usepackage{rank1}

\usepackage{color}
\usepackage{graphicx}
\input thmdefs.sty

\usepackage{mathrsfs}

\author{A.A.\,Prikhod'ko}
\address{Lomonosov Moscow State University}
\email{sasha.prihodko@gmail.com}

\def\muCL{\mu_{\mathrm{CL}}}
\def\muCLh{\Hat\mu_{\mathrm{CL}}}
\def\Rajchman{{\mathscr{R}}}

\def\BorelM{{\mathscr{M}}}
\def\Fourier{{\mathsf F}}
\def\wg{{\boldsymbol \gamma}}

\begin{document}

\title[Singular Schaeffer--Salem measures]
	{Singular Schaeffer--Salem measures on $[0,1]$ of~dynamical~system origin}
	
\begin{abstract} 
We study a class of dynamical systems given by measure preserving actions of 
the group $\Set{Z}^d$ or $\Set{R}^d$ and generating a set of spectral measures~$\sigma$ 
with an extremal rate of the Fourier coefficient decay: 
${ \Hat\sigma(n) = O(|n|^{-1/2+\eps}) }$ for any ${\eps > 0}$. 
Singular measures with this property are investigated in works due to 
Wiener and Wintner, Schaeffer, Salem, Ivashev-Musatov, Zygmund et al. 
Thus, the discovered effect provides 
a~new construction of singular distributions of Schaeffer--Salem type on the torus~$\Set{T}^d$ 
and in the space~$\Set{R}^d$. 
 
This work is supported by CNRS (France, Normandie), RFFI grant \No\,11-01-00759-a 
and the grant ``Leading Russian scientific schools'' \No\,NSh-5998.2012.1. 
\end{abstract}
 
\maketitle
  
\section[Constructions and properties of Borel measures]%
	{Constructions and properties of Borel measures on $[0,1]$}
 
In this work we investigate a class of dynamical systems generating 
spectral measures characterized by fast Fourier coefficient decay. 
We start with a well-known construction due to Riesz~\cite{RieszOnBinProd}.
He proposed to consider a~formal infinite product  
\begin{equation}\label{eRClassique} 
	\prod_{n=1}^\infty \bigl( 1 + a_n \cos(\omega_n x+\phi_n) \bigr), 
\end{equation} 
where $\omega_n \in 2\pi\Set{Z}$ is an increasing  sequence, 
${0 < a_n \le 1}$ and ${\phi_n \in \Set{R}}$. 
It is well known that for a certain choice of parameters $a_n$, $\omega_n$ and~$\phi_n$, 
for example, if ${\omega_{n+1}/\omega_n \ge q > 3}$ and ${\sum_n a_n^2 = \infty}$, 
this product represents a singular measure on~$[0,1]$ (see~\cite{ZygmundTSII},~\S\,7). 
We understand this statement as follows. The finite products 
$$
	\rho_N(x) = \prod_{n \le N} \bigl( 1 + a_n \cos(\omega_n x+\phi_n) \bigr) 
$$
are interpreted as densities of probability measures on $[0,1]$, and we have convergence 
$$
	\rho_N(x)\,ds \to d\sigma \quad \text{при} \quad n \to \infty 
$$ 
in the weak topology, where $\sigma$ is a measure on~$[0,1]$.   
The infinite products~\eqref{eRClassique} 
today refered to as {\it classical Riesz products\/}   
as well as generalized Riesz products 
$$
	\prod_{n=1}^\infty P_n(z), \qquad 
	P_n(z) = \sum_{k=0}^{q_n-1} c_{n,k}\, z^k, \quad z \in \Set{C}, \quad |z| = 1,  
$$
provide an important construction of singular measures 
broadly applied in analysis and dynamical systems (see, \cite{IsmagilovOnRieszPAndRandWalkRE}, 
	\cite{Bourgain}, \cite{ChoksiNadkarni}, 	
	\cite{Guenais}, \cite{Schaeffer}). 
 
Let us denote $\BorelM([0,1])$ or simply~$\BorelM$ 
the class of all Borel probability measures on~$[0,1]$. 
A~measure $\nu$ is {\it absolutely continuous\/} with respect to another measure~$\mu$, 
notation: ${\nu \ll \mu}$, 
if ${\nu = p(x)\mu}$, where $p(x)$ is a certain density, ${p(x) \in L^1(\mu)}$. 
The relation ${\nu \ll \mu}$ is a~partial order on~$\BorelM$. 
The measures $\mu$~and~$\nu$ in the class~$\BorelM$ are called {\it mutual singular}, ${\mu \perp \nu}$, 
if there exists a Borel set $E$ such that ${\mu(E) = \nu(E^c) = 1}$, 
where $E^c$ is the complement  to the set~$E$. 
Recall that any Borel measure $\sigma$ on a segment in the real line 
is uniquely expanded in a~sum 
\begin{equation}\label{eExpansionOfSigma}  
	\sigma = \sigma_d + \sigma_s + \sigma_{ac}, \qquad 
	\sigma_s \perp \la, \quad 
	\sigma_{ac} = p(x)\la, 
\end{equation} 
of discrete (a sum of atoms), purely singular and absolutely continuous components, 
where $\la$ is the normalized Lebesgue measure on the segment. 
Since every measure $\sigma \in \BorelM([0,1))$ can be considered as a measure on~$\Set{R}$ 
one can define the {\it Fourier transform\/} of~$\sigma$ by the formula
$$
	\Fourier[\sigma](t) = \int_0^{2\pi} e^{-2\pi i\, tx} \,d\sigma(x), \qquad t \in \Set{R}. 
$$
For probability measures on~$\Set{R}$ 
the following general observation holds  (see~\cite{ShiryaevE}, ch.\,2,~\S\,12). 
If the Fourier transform of~$\sigma$ belongs to $L^1(\Set{R})$, 
for example, if ${\Fourier[\sigma](t) = O(t^{-1-\alpha})}$, ${\alpha > 0}$, 
then $\sigma$ in absolutely continuous. At the same time, 
any measure $\sigma \in \BorelM([0,1))$ as a~measure 
on the compact group ${\Set{T} = \Set{R}/\Set{Z} \simeq [0,1)}$ 
generates the sequence of {\it Fourier coefficients\/} $\Hat\sigma(n)$ 
supported on the dual group ${\Hat{\Set{T}} = \Set{Z}}$.
Note that 
$$
	\Hat \sigma(n) = \Fourier[\sigma](n) = \int_0^{2\pi} e^{-2\pi i\, nx} \,d\sigma(x), 
	\qquad n \in \Set{Z}. 
$$
Suppose that ${\Hat\sigma \in l^2(\Set{Z})}$. Then the Fourier series 
$\sum_n \Hat\sigma(n)\,e^{2\pi i\,nt}$ converges in the space $L^2(\Set{T})$ to some 
function $p(x)$. Using Cauchy--Schwarz inequality 
$\|p\|_1 = \scpr<|p|,1> \le \|p\|_2$ we see that $p(x)$ is a density of some measure 
${p(x)\,dx \in \BorelM(\Set{R})}$,
and ${\Hat p(n) = \Hat\sigma(n)}$ (eg.\ see~\cite{WienerWintner}, \S\,1). 
Furhter, notice that any sequence ${ c_n = O(n^{-1/2-\alpha}) }$, ${\alpha > 0}$, 
is square summable. 
 
In the case of singular measure $\sigma$, as a~rule, we deal with 
a divergent series $\sum_n \Hat\sigma(n)\, z^n$. Understanding analytic properties 
of a~singular measure~$\sigma$, when we know certain com\-bina\-torial properties of the sequence $\Hat\sigma(n)$,   
becomes a very complicated problem. 
Louzin \cite{LuzinOnDivSeriesE} constructed the first example of 
power series $\sum_n c_n\, z^n$ with ${c_n \to 0}$ divergent everywhere on the unit circle ${|z| = 1}$. 
Further, Neder \cite{NederOnDivSer} proved that any series $\sum_n c_n\, z^n$ 
with the property ${\sum_n |c_n|^2 = \infty}$ can be transformed to everywhere divergent 
(for ${|z| = 1}$) using some phase correction $\tilde c_n = e^{i\,\phi_n}\,c_n$. 
Now let us turn to expansion~\eqref{eExpansionOfSigma} of~$\sigma$ and remark that 
Riemann--Lebesgue lemma can be interperted in the folowing way. 

\begin{lem}
Given an absolutely continuous measure $\sigma_{ac}$, 
$$\Hat\sigma_{ac}(n) \to 0 \quad \text{as} \quad n \to \infty.$$ 
\end{lem}
 
\begin{defn}
We call {\it Menshov--Rajchman measure}, a singular measure satisfying
$\Hat\mu(n) \to 0$, ${n \to \infty}$. 
We denote as $\Rajchman$ the class of all measures of such kind. 
\end{defn}

Evidently any discrete measure $\sigma_d$ is never of Menshov--Rajchman type since 
its Fourier transform $\Hat\sigma_d(n)$ is a Bohr almost periodic sequence. 
At the same time, it is easy to see that the singular Cantor--Lebesgue measure $\muCL$ 
supported on the standard $1/3$-Cantor set enjoys the property ${\muCLh(3^k n) = \muCLh(n)}$, 
which is explained by fractal symmetry of the Cantor set. Thus, ${\muCL \not\in \Rajchman}$. 
Modifying the construction of~$\muCL$ Menshov \cite{Menshov} 
provided the first example of singular measure in the class~$\Rajchman$. 
Further, Neder \cite{NederOnRajchman} proved that any Menshov--Rajchman measure 
cannot be a mixture of discrete and continuous component, and then Wiener \cite{WienerOnQVar} 
extended this result and showed that the Fourier coefficients of any continuous measure~$\mu$ 
converge to zero in average, and any set $\{n \where |\Hat\mu(n)| > b > 0 \}$ 
has zero density in~$\Set{Z}$. 
Littlewood \cite{LittlewoodOnPowRajchman} found 
a singular probability measure~$\sigma$ 
with the rate of decay 
$$
	\Hat\sigma(n) = O\bigl( |n|^{-c} \bigr), \qquad c > 0. 
$$
Then Wiener and Wintner \cite{WienerWintner} 
obtained a stronger result demonstrating 
that the exponent $c$ can be arbitrary close to~$1/2$, 
but the approach proposed by the authors 
generates a measure $\sigma$ that depends on ${c = 1/2-\alpha}$, $\alpha > 0$. 
Soon after this work Schaeffer~\cite{Schaeffer} 
using the idea of Riesz products 
proved the existence of a singular $\sigma$ with 
$$
	\Hat\sigma(n) = O\bigl( r(|n|)\cdot |n|^{-1/2} \bigr) 
$$ 
for any given increasing sequence ${r(n) \to \infty}$, ${n \to \infty}$. 
In particular, $\sigma$ 
satisies 
$$
	{\Hat\sigma(n) = O(|n|^{-1/2+\eps})} \quad \text{for any} \quad \eps > 0. 
$$
Ivashev-Musatov \cite{IvasevMusatovOnSingF} got a further improvement of Schaeffer's result. 
He found a set of singular measures with sub-$|n|^{-1/2}$ rate of correlation decay satisfying 
${\Hat\sigma(n) = O(\rho(n)\cdot |n|^{-1/2+\eps})}$ with ${\rho(n) \to 0}$ 
but ${\rho(n) \gg |n|^{-\eps}}$ for any ${\eps > 0}$. 
Following~\cite{WienerWintner} let us denote $\kappa(\sigma)$ the infinum of 
real~$\gamma$'s such that ${\Hat\sigma(n) = O(|n|^\gamma)}$.
In a~series of works \cite{SalemOnStrMonSingF,SalemOnSetsUniqueness,SalemOnSetsUniquenessII} 
Salem introduced an approach that helps to see, in particular, explicit examples 
of singular distributions with the property ${\kappa(\sigma) = -1/2}$. 

\begin{defn}
Let us call singular measures on $\Set{T}^d$ (respectively,  $\Set{R}^d$)  
satisfying the condition $\Hat\sigma(n) = O(|n|^{-d/2+\eps})$ for any ${\eps > 0}$, 
{\it measures of Schaeffer--Salem type}. 
\end{defn}
 
In this note we discover that Schaeffer--Salem measures appear as spectral measures 
for a class of group actions with invariant measure. 
Let us consider a measure preserving invertible transformation $T \Maps X \to X$ 
of the standard Lebesgue space $(X,\cB,\mu)$ and define {\it Koopman operator\/} 
on the space ${H = L^2(X,\cB,\mu)}$, 
$$
	\hat T \Maps H \to H \Maps f(x) \mapsto f(Tx). 
$$
Clearly, $\hat T$ is a unitary operator in 
a separable Hilbert space $H$, hence, it is characterized up to a unitary equivalence 
by the pair $(\sigma_{(T)},\Mult(z))$, where $\sigma_{(T)}$ is the 
{\it measure of maximal spectral type\/} and 
$\Mult(z)$ is the {\it multiplicity function}. Of course, two transformations 
which are spectrally isomorphic, need not be isomorphic as dynamical systems. 
For example, all Bernoulli shifts	 have Lebesgue spectrum of infinite multiplicity 
but they are dinstiguished by entropy. 
And, in~fact, it is a hard problem far from complete understanding 
to classify all pairs $(\sigma_{(T)},\Mult(z))$ that can appear as spectral invariants 
of a~measure preserving transformation (see~\cite{KatokThouvenotSpTh,LemEncycloSpTh}). 
Further, given an element $f \in L^2(X,\cB,\mu)$, the {\it spectral measure\/} 
${\sigma_f \in \BorelM}$ on $\Set{T} \simeq [0,1)$ is uniquely defined by the relation 
$$
	\Hat\sigma_f(n) = \int e^{-2\pi i\, xn} \,d\sigma(x) = R_f(-n) \stackrel{def}{=} \scpr<\hat T^{-n}f,f>. 
$$
It is easy to see that $\sigma_f \ll \sigma_{(T)}$.

\section[Schaeffer--Salem spectral measures]
	{Dynamical systems generating Schaeffer--Salem measures} 

In this section we define actions of the groups $\Set{Z}^d$ and $\Set{R}^d$ with invariant measure 
generating spectral measures of Schaeffer--Salem type for a dense set of function on the phase space. 
Without loss of generality we concentrate our attention on the case of $\Set{Z}^d$-actions. 
Consider a~nested sequence of lattices $\G_n$ in the acting group ${G = \Set{Z}^d}$ 
having ${G_{n+1} \subset G_n}$, and let ${M_n = G/\G_n}$ be the corresponding sequence 
of homogeneous spaces linked with a~natural projection $\pi_n \Maps M_{n+1} \to M_n$ 
mapping $a + \G_{n+1}$ onto the point $a + \G_n$. 
Let us also fix a F\"olner sequence of $\G_n$-fundamental domains $U_n$, where 
$$
	G = \bigcup_{\gamma \in \G_n} (\gamma + U_n)
\quad \; \text{and} \; \quad 
	\gamma_1 + U_n \cap \gamma_2 + U_n = \eset, 
	\quad \text{whenever} \quad \gamma_1 \not= \gamma_2. 
$$
To simplify understanding of the construction let us consider 
a particular case $\G_n = h_n\Set{Z}^d$, ${h_{n+1} = q_n h_n}$, $q_n \in \Set{Z}$, 
and let $U_n$ be rectangles ${U_n = [0,h_n)^{\times d}}$. 
Now let us introduce the maps 
$\phi_n \Maps M_{n+1} \to M_n$, 
$$
	\phi_n(\gamma + u) = u + \a_{n,\gamma}, \qquad \gamma \in \G_n / \G_{n+1}, \quad u \in U_n. 
$$
The idea of $\phi_n$ can be explained as fllows. Take a pair $M_n$ and $M_{n+1}$. 
We consider $M_{n+1}$ as a~union (up to a~boundary in the case of~$\Set{R}^d$) 
of finitely many domains $\gamma + U_n$, where ${\gamma \in \G_n / \G_{n+1}}$. 
Next we note that 
each domain $\gamma + U_n$ is mapped by~$\pi_n$ in $1$-to-$1$ way 
to the manifold $M_n$, and then we define $\phi_n$ in such a way that 
$\gamma + U_n$ after descending to~$M_n$ is ``rotated'' to some fixed vector $\a_{n,\gamma}$. 
In a sense the boundaries of the domains $\gamma + U_n$ contains all 
the discontinunity points of the map $\phi_n$. 
For a given point ${a \in M_{n+1}}$ let us define $\wg_n(a)$ to be the unique element 
$\gamma \in \G_n / \G_{n+1}$ such that ${ a \in \gamma + U_n}$. 
Let $\mu_n$ be the normalized shift invariant measure on the space $M_n$. 
Then $\mu_{n+1}(\phi_n^{-1}(A)) \equiv \mu_{n}(A)$. Further, let us define $X$ 
to be the inverse limit of the spaces~$M_n$, 
$$
	X = \bigl\{ x = (x_1,x_2,\dots,x_n,\dots) \where x_n \in M_n,\ \phi_n(x_{n+1}) = x_n \bigr\}. 
$$
The space $X$ is endowed both with the weak topology and with the structure of probability space 
$(X,\cB,\mu)$, where $\sigma$-algebra $\cB$ is generated by the cylindric sets 
$B_{n,A} = \{x \where x_n \in A\}$, where $A$ is a Borel set in~$M_n$, 
and $\mu$ is the measure extending $\mu_n$. Notice that ${\tilde\pi_n^* \mu = \mu_n}$, 
where ${\tilde\pi_n(x) = x_n}$. 

Finally let us define an action of the group $G$ on the space $(X,\cB,\mu)$. 
Assume now that ${q_{n+1} > 2q_n}$. 
Fix $t \in G$. It is to check using Borel--Cantelli lemma that the probability 
$$
	\mu \bigl\{ 
		x \in X \where 
		\exists\,n_0(x)\quad \forall\,n \ge n_0(x) \quad \wg_n(x_n) = \wg_n(t+x_n)      
	\bigr\} 
	= 1, 
$$
hence, for the points $x$ of such kind (belonging to the set above) 
we can define $T^t$ by the rule $(T^tx)_n = t+x_n$ for ${n \ge n_0(x)}$. 
For indexes ${n < n_0(x)}$ the coordinates $(T^tx)_n$ are recovered 
using the fundamental equation $\phi_n(x_{n+1}) = x_n$. 
The following lemma directly follows from the definition 
(see also \cite{IcePaperI} for the careful examination of the case ${G = \Set{Z}}$). 

\begin{lem}
The maps $T^t$ provide an action of the group $G$ 
on the space $(X,\cB,\mu)$ by measure preserving transformations. 
\end{lem}
 
We call the class of systems provided by this construction {\it systems of iceberg type}. 
It can be easilily seen that this class is analogous to the general rank one actions 
(see \cite{DanilenkoOnFunnyRankOneWMix,Prikhodko}). The spectral measures 
of rank one transformation and flows are given by a certain class of generalized Riesz products 
\cite{Bourgain,ChoksiNadkarni,LebesgueFlows} and for the study of spectral type 
for systems of iceberg type one can also use a variation of Riesz products.

Let us remark that the $\Set{Z}$-action of iceberg type with $q_n = 2$, ${(\a_{n,0},\a_{n,1}) = (0,h_n/2)}$, 
is identical to the classical Morse transformation (see~\cite{VershikOnMorseArithm} for 
the discussion of arithmetic properties of Morse systems). 

\begin{thm}\label{thSalemSpMeasures} 
Let $\a_{n,\gamma}$ be a family of independent random variables, 
uniformely distributed on finite sets $\G_n / \G_{n+1}$. 
Then for a certain sequence $q_n \to \infty$ and for a set of (cylindric) functions $f$ 
dense in $L^2(X,\mu)$ the spectral measures $\sigma_f$ 
satisfy almost surely the condition 
\begin{equation}\label{eSalemSpMeasures} 
	\Hat\sigma_f(t) = O(|t|^{-1/2+\eps})
\end{equation}
for any $\eps > 0$. Moreover, $\Hat\sigma_f \not\in L^2(G)$, and 
$\Hat\sigma_f \in L^p(G)$ for ${p > 2}$. 
\end{thm}

In the case ${G = \Set{Z}}$ the proof of the next theorem is given in~\cite{IcePaperI}. 
 
The most important effect, underlying the proof of the theorem, 
is the universal character of the rate of correlation decay 
$R_f(-t) = \scpr<T^{-t}f,f> = \Hat\sigma_f(t)$ 
for this class of dynamic systems. 
As we will see, it follows easily from the calculation 
of expectations $\E |R_f(-t)|^2$ for random variables $R_f(-t)$, 
but reconstruction of $\sigma_f$ from the sigma $\Hat\sigma_f(t)$ is a complicated problem, 
and even for the random family of dynamical systems of theorem~\ref{thSalemSpMeasures} 
we cannot say whether measures $\sigma_f$ are almost surely singular? 

In paper \cite{VRyzhOnInfRankOneMix} Ryzhikov proves the existence of transformations 
preserving an infinite measure satisfying the estimate $\kappa(\sigma_f) \le -1/2$ 
for a dense set of elements $f$. 

\medskip
{\it Proof of theorem~\ref{thSalemSpMeasures}.} 
Without loss of generality let us assume that $d=1$. 
As it is mentioned above the proof of theorem~\ref{thSalemSpMeasures} 
is rather simple. 
It~is based on an elementary and common phenomenon which hides 
more special combinatorial properties of the correlation sequence $R_f(t)$ 
as well as properties of the spectral measure~$\sigma_f$. 
Indeed, a~universal estimate of the expectation $\E |R_f(t)|^2$ is established for our class of systems: 
$$ 
	\E |R_f(t)|^2 \le r(n)\cdot t^{-1}
$$ 
for 
${h_{n-1} \ll t \le h_n}$, where $r(n)$ is a steadily increasing sequence: ${ r(n) \ll h_n }$. 
Let $f$ be a cylindric function depending only on $x_{n_0}$, i.e.\ ${f(x) = f_{n_0}(x_{n_0})}$, 
It is enough to consider just the {\it correlation functions\/} 
$$ 
	R^\circ_n(t) = \frac1{h_n^d} \int_{M_n} f_n(x-t)\,\OL{f_n(x)} \,d\mu_n, 
$$
where $f_n$ is the lift of the function~$f$ to the manifold~$M_n$. 
It is also sufficient to consider a particular set of values 
${t = h_{n-1} \cdot s \in \G_{n-1}}$, ${s \in \Set{Z}^d \sms \{0\} }$. 
Such values of~$t$ are characterized by the following properties. 
Whenever $t \in \G_{n-1}$ the shift transformation $x \mapsto t + x$ 
preserves the partition of the manifold~$M_n$ to the sets $\gamma + U_{n-1}$. 
Then given a cylindric function $f$ with zero mean, $\int f\,d\mu = 0$, we see that 
$$
	\E R^\circ_n(t) = 
	\frac1{q_{n-1}} \sum_{\gamma \in \G_{n-1} / \G_{n}} 
		\E\scpr<\rho_{\alpha_{n-1,\gamma}} f_{n-1}, \rho_{\alpha_{n-1,\gamma+s}} f_{n-1}> = 0, 
$$
since $\alpha_{n-1,\gamma}$ are $\alpha_{n-1,\gamma+s}$ are independent. 
Developping this technique we establish a reccurent property 
$$
	\E |R^\circ_n(t)|^2 \sim \frac1{q_{n-1}} \cdot \E_{\mu_{n-1}} \E |R^\circ_{n-1}|^2, 
$$ 
and the estimate 
$$
	\E\|R^\circ_n\|^2 \le \E\|R^\circ_{n-1}\|^2 + h_n \cdot h_{n-1}^{-1}\E \|R^\circ_{n-1}\|^2  
	\le 2\:\E \|R^\circ_{n-1}\|^2 \cdot (1+o(1)) 
	= O(2^n), 
$$
hence, 
$$
	\E |R^\circ_n(t)|^2 = O(q_{n-1}^{-1} \cdot 2^n h_{n-1}^{-1}) 
	= O(h_n^{-1} \cdot 2^n ) = O(t^{-1} \cdot 2^n), 
$$
whereе $t \le h_n$ and the proof is finished. $\Box$

\medskip
It follows directly from theorem~\ref{thSalemSpMeasures} that 
${ \sigma_f \conv \sigma_f \ll \la }$, where $\la$ is the normalized invariant measure on~$\Hat G$. 
Thus, our observation is connected to the open question due to Banach, --- 
{\it ``Does there exist a $\Set{Z}$-action with invariant probability measure 
	having Lebesgue spectrum of multiplicity one?''} --- 
since the spectral multiplicity for almost every action in theorem~\ref{thSalemSpMeasures} equals~one 
(and for {\em all\/} systems of iceberg type we have ${\Mult(z) \le 4}$).

\medskip
{\it Open questions and hypotheses.} 
\begin{itemize} 
	\item[(i)] Dynamical systems in theorem~\ref{thSalemSpMeasures} 
		are almost surely of singular spectrum 
	\item[(ii)] $\kappa(\sigma_f) = -1/2$ for all systems of iceberg type, 
		i.e.\ for any choice of the parameters $\a_{n,\gamma}$, where $f$ is a cylindric function 
	\item[(iii)] Can we find a speed of Fourier coefficient decay $\Hat\sigma(t) = O(|t|^{-1/2}\,\rho(t))$, 
		that can be reached in the class of {\it all\/} singular measures, but {\em impossible\/} 
		for measures of maximal spectral type generated by measure preserving transformations? 
		(cf.\ question~(iv)) 
	\item[(iv)] Given a measure of Salem type, $\kappa(\sigma) = -1/2$, is it posible 
		to reach any speed of decay of type $\Hat\sigma(t) = O(|t|^{-1/2}\,\rho(t))$ 
		with ${ \forall\,\eps > 0\quad |t|^{-\eps} \le \rho(t) \le |t|^\eps }$ 
		just multiplying by some density ${p(x) \in L^1(\sigma) \:}$? 
	\item[(v)] Is it true that all $\Set{Z}^d$- and $\Set{R}^d$-actions of iceberg type have singular spectrum? 
\end{itemize} 

\medskip
The author is grateful to D.V.\,Anosov, A.M.\,St\"epin, A.M.\,Vershik, V.V.\,Ryzhikov, 
J.-P.\,Thouvenot, M.\,Lemanczyk, 
B.M.\,Gurevich, V.I.\,Oseledec, S.A.\,Pirogov, O.N.\,Ageev, I.A.\,Bogaevskii, 
S.V.\,Konyagin, I.D.\,Shkredov, G.\,Alkauskas и A.V.\,Egorov for fruitful discus\-sions and helpful remarks.

\bibliographystyle{amsplain}
\bibliography{IcePaperI}
 
\end{document}